\documentclass[11pt]{article}
\usepackage{amsmath,amssymb,amsfonts,amsthm,latexsym,amstext,amscd}
\usepackage[english]{babel}
\usepackage{fancyhdr}

\newtheorem{theorem}{Theorem}

\newtheorem*{conj}{Conjecture}
\newtheorem*{theoA}{Theorem A}

\theoremstyle{definition}

\newtheorem{remark}{Remark}

\newtheorem*{rema}{Remark}

\newcommand{\beq}{\begin{equation}}
\newcommand{\eeq}{\end{equation}}

\numberwithin{equation}{section}
\begin{document}

\title{On a conjecture of Erd\H{o}s, P\'olya and Tur\'an on consecutive gaps between primes}
\author{J\'anos Pintz\thanks{Supported by ERC-AdG. 321104 and OTKA NK104183.}}


\date{}



\maketitle

\section{Introduction}
\label{sec:1}

One of the favourite problems of Erd\H{o}s (and Tur\'an) was to investigate local problems in the distribution of primes, in particular to examine gaps or blocks of successive gaps between consecutive primes.

Let $\mathcal P := \{p_n\}_1^\infty$ be the sequence of all primes and
\beq
\label{eq:1.1}
d_n = p_n - p_{n - 1} \ \ \ (n = 2, 3, \dots)
\eeq
be the sequence of gaps between consecutive primes.

In 1948 Erd\H{o}s and Tur\'an \cite{ET} showed that
\beq
\label{eq:1.2}
d_{n + 1} - d_n
\eeq
changes sign infinitely often.
Soon after this Erd\H{o}s \cite{Erd1} showed the stronger relation
\beq
\label{eq:1.3}
\liminf\limits_{n \to \infty} \frac{d_{n + 1}}{d_n} < 1 < \limsup_{n \to \infty} \frac{d_{n + 1}}{d_n} .
\eeq

In the same work \cite{ET}, that is, already 67 years ago, Erd\H{o}s and Tur\'an asked for a necessary and sufficient condition that
\beq
\label{eq:1.4}
\sum_{i = 1}^k a_i p_{n + i}
\eeq
should have infinitely many sign changes as $n \to \infty$, where $a_1, \dots, a_k$ are given real numbers.
They observed that
\beq
\label{eq:1.5}
\sum_{i = 1}^k a_i = 0
\eeq
is clearly necessary, and P\'olya observed that if \eqref{eq:1.4} has infinitely many sign changes, then the $k$ numbers
\beq
\label{eq:1.6}
\alpha_j = \sum_{i = 1}^j a_n
\eeq
cannot all have the same sign.
As described in \cite{ET} and \cite{Erd2}, Erd\H{o}s, P\'olya and Tur\'an then conjectured that the above condition on $\alpha_j$ is a necessary and sufficient condition for the infinitely many sign changes of \eqref{eq:1.4}.
As Erd\H{o}s writes on p.\ 12 of \cite{Erd2}:
``We are very far from being able to prove this, in fact I cannot even prove that $d_n > d_{n + 1} + d_{n + 2}$ has infinitely many solutions.
I proved the following much easier theorem:
Assume that
\beq
\label{eq:1.7}
\sum_{i = 1}^{k - 1} \alpha_i = 0 \ \ \text{ and } \ \  \alpha_{k - 1} \neq 0.
\eeq
Then \eqref{eq:1.4} changes sign infinitely often.''

In my recent work \cite{Pin} I showed several partial results in this direction (see Theorems 17--19) but I was far from being able to show the original conjecture of Erd\H{o}s, P\'olya and Tur\'an.
In the present work I will show the original conjecture based on the recent groundbreaking ideas of J. Maynard \cite{May} and T. Tao \cite{Tao} on bounded gaps between primes.

\section{Some remarks and a stronger form of the Erd\H{o}s--P\'olya--Tur\'an conjecture}
\label{sec:2}

Since the necessity of \eqref{eq:1.5} is trivial we can further on always suppose \eqref{eq:1.5}.
So we can rewrite \eqref{eq:1.4} as
\begin{align}
T &= \sum_{i = 1}^k a_i p_{n + i} = \sum_{i = 1}^k a_i \biggl(p_n + \sum_{\nu = 1}^i d_{n + \nu}\biggr) = \sum_{i = 1}^k a_i \biggl( \sum_{\nu = 1}^i d_{n + \nu}\biggr)\nonumber \\
&= \sum_{j = 1}^k d_{n + j} \sum_{i = j}^k a_i = -\sum_{j = 1}^k d_{n + j} \alpha_{j - 1} = - \sum_{j = 2}^k d_{n + j} \alpha_{j - 1}
\label{eq:2.1}
\end{align}
if we define $\alpha_0 = 0$.

Thus the original conjecture is equivalent to the following one (if we let $\ell := k - 1$).

\begin{conj}[Erd\H{o}s--P\'olya--Tur\'an]
The expression
\beq
\label{eq:2.2}
\sum_{i = 1}^\ell \alpha_i d_{n + i}
\eeq
changes sign infinitely often as $n$ runs through all integers if and only if the non-zero elements among $\alpha_1, \alpha_2, \dots, \alpha_\ell$ do not all have the same sign.
\end{conj}

\begin{theorem}
\label{th:1}
The above conjecture is true.
\end{theorem}

The above theorem clearly follows from (but as it is easy to see, is in fact equivalent to)
the following one.

\begin{theorem}
\label{th:2}
We have for every fixed natural number $\ell$
\beq
\label{eq:2.3}
\limsup_{m \to \infty} \frac{d_m}{\max\bigl(d_{m - \ell}, \dots, d_{m - 1}, d_{m + 1}, \dots , d_{m + \ell}\bigr)} = \infty .
\eeq
\end{theorem}

We will prove this in the stronger form expressed by

\begin{theorem}
\label{th:3}
For every natural number $\ell$ there exists an explicitly calculable constant $c(\ell) > 0$ such that
\beq
\label{eq:2.4}
\limsup_{m \to \infty} \frac{d_m}{(\log m)^{c(\ell)} \max\bigl(d_{m - \ell}, \dots, d_{m - 1}, d_{m + 1}, \dots, d_{m + \ell}\bigr)} > 0.
\eeq
\end{theorem}

\begin{rema}
It follows from the proof that one can take
\beq
\label{eq:2.5}
c(\ell) = C_1 e^{-C_2 \ell}
\eeq
with positive absolute constants $C_1, C_2 > 0$.
\end{rema}

\section{Proof of Theorem~\ref{th:3}}
\label{sec:3}

The proof will be based on the first assertion (see (i)) of a very nice result of W. D. Banks, T. Freiberg and J. Maynard \cite{BFM} which appears as Theorem~4.3 in their work, which we quote now restricted on (i) and with a slight change as

\begin{theoA}[Banks--Freiberg--Maynard]
 Let $m$, $k$ and $\epsilon = \epsilon(k)$ be fixed.
 If $k$ is a sufficiently large multiple of $16m + 1$ and $\epsilon$ is sufficiently small, there is some $N(m, k, \epsilon)$ such that the following holds for all $N \geqslant N(m, k, \epsilon)$.
 With $Z_{N^{4\epsilon}}$ given by {\rm (4.8)} of \cite{BFM}, let
 \beq
 \label{eq:3.1}
 w = \epsilon \log N \ \text{ and } \ W = \prod\limits_{\substack{p\leqslant w\\ p \nmid Z_{N^{4\epsilon}}}} p.
 \eeq
 Let $\mathcal H = \{h_1, \dots, h_k\}$ be an admissible $k$-tuple such that
 \beq
 \label{eq:3.2}
 0 \leqslant h_1, \dots, h_k \leqslant N
 \eeq
 and
 \beq
 \label{eq:3.3}
 p \Big| \prod\limits_{1 \leqslant i < j \leqslant k} (h_j - h_i) \Longrightarrow p \leqslant w.
 \eeq
 Let
 \beq
 \label{eq:3.4}
 \mathcal H = \mathcal H_1 \cup \cdots \cup \mathcal H_{16 m + 1}
 \eeq
 be a partition of $\mathcal H$ into $16 m + 1$ sets of equal size.
 Finally, let $b$ be an integer such that
 \beq
 \label{eq:3.5}
 \biggl(\prod\limits_{i = 1}^k (b + h_i), W \biggr) = 1.
 \eeq
 There is some $n_1 \in (N, 2N]$ with $n_1 \equiv b \text{ {\rm mod} } W$, and some set of $m + 1$ distinct indices $\{i_1, \dots, i_{m + 1}\} \subseteq \{1, \dots, 16m + 1\}$, such that
 \beq
 \label{eq:3.6}
\bigl|\mathcal H_i(n_1) \cap \mathbb P \bigr| = 1 \ \text{ for all } \ i \in \bigl\{i_1, \dots, i_{m + 1}\bigr\}.
 \eeq
\end{theoA}

\begin{remark}
\label{rem:1}
The definition of $Z_{N^{4\epsilon}}$ is given earlier in the work \cite{BFM} but its value does not play a significant role in the application of the result (it is the greatest prime factor of a possible exceptional modulus if such a modulus exists and it is equal to $1$ if no such modulus exists).
\end{remark}

\begin{remark}
\label{rem:2}
According to the calculation of the present author $8m + 1$ in (4.18) of \cite{BFM} has to be replaced by $16m + 1$.
\end{remark}

The proof uses the Maynard--Tao method \cite{May}, \cite{Tao} and other important ideas as a modified Erd\H{o}s--Rankin type construction (see Section~5 of \cite{BFM}), a modified Bombieri--Vinogradov theorem, somewhat similar to Theorem~6 of \cite{GPY}, and an important observation of the Polymath project \cite{Pol} according to which one can estimate from above how often we have more than one prime in the translation of a subset $\mathcal H'$ of an admissible $k$-tuple $\mathcal H$ (in the weighted sense).

We note that the variable $k$ in Theorem~A has nothing to do with the one appearing in Section~1 and \eqref{eq:2.1} of our work which satisfied $k = \ell + 1$.
The present $k$ will be here a large multiple of $16m + 1$ and $m$ will satisfy $m \asymp \ell$ here.
In fact we will define now
\beq
\label{eq:3.7}
L := \ell + 2, \ \ \ m := 62L - 33 = 62 \ell + 91.
\eeq
The Maynard--Tao method needs to choose in the proof of Theorem~A
\beq
\label{eq:3.8}
k = \exp (C_3 m) = \exp(C_4 \ell)
\eeq
(by the relation $\delta \varrho \log k = 2m$, appearing in the first line on p.\ 17 of \cite{BFM}).
This will imply the appearance of $(\log n)^{c(\ell)}$ in \eqref{eq:2.4} of us with $c(\ell)$ defined as in \eqref{eq:2.5}, that is $C_1 \exp(-C_2 \ell)$.

In order to show Theorem~\ref{th:3} we will choose with a sufficiently large $k$ an admissible $k$-tuple of $\mathcal H$ with
\beq
\label{eq:3.9}
2(16m + 1) \mid k,
\eeq
for every given sufficiently large~$N$.
We further let
\beq
\label{eq:3.10}
J := 32L - 17
\eeq
which implies
\beq
\label{eq:3.11}
16m + 1 = 992L - 527 = 31J .
\eeq

We will partition our admissible $k$-tuple $\mathcal H$ into $16m + 1 = 31J$ subsets of equal size $k/(31J)$.

We will use the additional information of \cite{BFM} (see Sections 5 and 6 of it) that by the Erd\H{o}s--Rankin procedure one can find for any sufficiently large $N$ an admissible $k$-tuple $\mathcal H$ and a number $n \in [N, 2N]$ which we fix in the following, such that with a
\beq
\label{eq:3.12}
z > \log N \sqrt{\log_2 N}
\eeq
all numbers of the form
\beq
\label{eq:3.13}
n + \nu, \ \ 1 < \nu \leqslant z, \ \ \nu \notin \mathcal H
\eeq
should be composite.
Hence all possible primes in $(n + 1, n + z]$ should be of type $n + h_i$, $h_i \in \mathcal H$.
We have here a lot of freedom in choosing $\mathcal H$.
First its elements can be as large as $\log N$ and the conditions of Theorem~A allow us to choose its  elements as
\beq
\label{eq:3.14}
h_i = \sum_{j = 1}^i b_j \ \ \text{ and } \ \ b_i = (1 + o(1)) \beta_j(\log N)^c j
\eeq
for any choice of $\beta_1,\dots, \beta_k$ and $c_1,\dots, c_k \in (0, 1]$
(see Sections 5-6 of \cite{BFM}).
We will choose $\beta_i = 1$ ($i = 1,2, \dots, k$) and with
\beq
\label{eq:3.15}
K = \frac{k}{62J}\,, \ \ i = (\nu J + \mu)K + \lambda, \ \ 0 \leqslant \nu \leqslant 30, \ \ 0 \leqslant \mu \leqslant J - 1, \ \ 1 \leqslant  \lambda \leqslant  K
\eeq
we will choose
\beq
\label{eq:3.16}
c_i = f(\nu, \mu, \lambda) := \frac{(30 - \nu) J + \mu}{62J} + \frac{\lambda}{k} \in \left[\frac1{k}, \frac12\right].
\eeq
This means that
\beq
\label{eq:3.17}
h_j - h_i = \sum_{t = i + 1}^j b_t \sim \max_{i < t \leqslant j} b_t \ \ \text{ for } \ i < j.
\eeq
From \eqref{eq:3.16} we further see that $c_i = f(\nu, \mu, \lambda)$, and consequently $b_i$ will be monotonically increasing in both $\mu$ and $\lambda$ for each fixed value of $\nu$ when $\mu \in [0, J - 1]$, $\lambda \in [1, K]$.
More exactly, for every fixed $\nu$ we have
\beq
\label{eq:3.18}
f(\nu, \mu, \lambda) - f(\nu, \mu', \lambda') \geqslant \frac1{k} \ \text{ if } \ \mu > \mu' \ \text{ or } \ \mu = \mu', \ \lambda > \lambda'.
\eeq
On the other hand this construction shows that
\beq
\label{eq:3.19}
f(\nu_2, \mu_2, \lambda_2) - f(\nu_1, \mu_1, \lambda_1) = \frac{\nu_1 - \nu_2}{62} + \frac{\mu_2 - \mu_1}{62J} + \frac{\lambda_2 - \lambda_1}{k} \leqslant - \frac1{k}
\eeq
if $\nu_2 > \nu_1$ for every quadruple $(\mu_1, \mu_2, \lambda_1, \lambda_2)$ if $\mu_i \in [0, J - 1]$, $\lambda_i \in [1, K]$
which means that $f(\nu, \mu, \lambda)$ is monotonically decreasing in $\nu$ independently of the values of $\mu$ and $\lambda$.

Let us define now the partition of $\mathcal H$ into $31J = 16m + 1$ subsets $\mathcal H_{\nu, \mu}$ $(0 \leqslant \nu \leqslant 30$, $0 \leqslant  \mu \leqslant  J - 1)$ as
\beq
\label{eq:3.20}
\mathcal H_{\nu, \mu} := \{h_i\}_{i \in I_{\nu, \mu}} \ \ I_{\nu, \mu} = \bigl\{i = (\nu J + \mu) K + \lambda, \ 1 \leqslant \lambda \leqslant K \bigr\}
\eeq
and let us organize these subsets into 31 columns according to the value of the index $\nu$ for $\nu = 0, 1, \dots, 30$.

The observations \eqref{eq:3.17}--\eqref{eq:3.19} show that the values of $b_i$ are increasing (by a factor $\geqslant (1 + o(1)) (\log N)^{1/k}$) within each column.
Further we see that if $b_i$ is in another column than $b_j$ with an index
\beq
\label{eq:3.21}
\nu(i) < \nu(j) ,
\eeq
then necessarily
\beq
\label{eq:3.22}
\frac{b_i}{b_j} \geqslant  (\log N)^{1/k}(1 + o(1)).
\eeq

This means by \eqref{eq:3.17} that if $n + h_i$ and $n + h_j$ $(i < j)$ are consecutive primes, then their difference is asymptotically equal to the dominant $b_t$ with $i < t \leqslant j$ and the ratio between two consecutive primegaps will be
\beq
\label{eq:3.23}
\geqslant (\log N)^{1/k}(1 + o(1)) \ \text{ or } \ \leqslant (1 + o(1)) (\log N)^{-1/k}.
\eeq

Theorem~\ref{th:2} will be shown if we can reach in one of the columns with index $\nu = 0,1, \dots, 29$ (that is, $\nu \neq 30$) at least $L$ primes of the form $n + h_i$ in such a way that we should have still in total at least $L$ primes of type $n + h_i$ in all remaining columns with an index larger than~$\nu$.

In this case we can choose the largest index $i$ within that column (that is, with $\nu(i) = \nu$) as our $h_i$ for which $n + h_i \in \mathcal P$, and we let
\beq
\label{eq:3.24}
p_m := n + h_i.
\eeq

This will imply that we have additionally $p_{m - 1}, \dots, p_{m - L + 1} = p_{m - \ell - 1}$ in the same column and the differences $d_m, d_{m - 1}, \dots, d_{m - \ell}$ satisfy
\beq
\label{eq:3.25}
\frac{d_{m - i}}{d_{m - j}} \gg (\log N)^{(j - i)/k} \ \text{ for } \ 0 \leqslant i < j \leqslant \ell
\eeq
in accordance with \eqref{eq:2.4}.

Further, in view of \eqref{eq:3.21}--\eqref{eq:3.22}, as all the later primes of type $p_{m + t}$ with $t \geqslant 1$ are of the form $n + h_j$ with $\nu(j) > \nu(i)$ we will have for the increments the relation \eqref{eq:3.22} and this will yield
\beq
\label{eq:3.26}
\frac{d_m}{d_{m + t}} \gg (\log N)^{1/k} \ \text{ for } \ 0 < t, \ p_{m + t} \leqslant n + z.
\eeq

So, let us suppose now that the first column having at least $L$ primes of the form $n + h_i$ has index $y$, where $0 \leqslant y \leqslant 30$.
If such an index, that is, such a column does not exist, then we have in total at most
\beq
\label{eq:3.27}
31L < 62 L - 32 = m + 1
\eeq
primes among $n + h_i$ in contradiction with \eqref{eq:3.6} in Theorem~A.
So we have such a column with index $y \in [0, 30]$.
This column contains at most $J$ subsets of type $\mathcal H_{i_t}$ described in \eqref{eq:3.6}.
If we have no further column at all (i.e.\ $y = 30$) or the number of primes in later columns is in total at most $L - 1$, then we have in total at most $30(L - 1)$ primes in all other columns.
This means that the total number of subsets $\mathcal H_{i_t}$ with exactly one prime of the form $n + h_j$ in it $(h_j \in \mathcal H_{i_t})$ is at most (cf.\ \eqref{eq:3.6}, \eqref{eq:3.7} and \eqref{eq:3.10})
\beq 
\label{eq:3.28}
30(L - 1) + J = 62L - 47 < 62L - 32 = m + 1
\eeq
in contradiction with \eqref{eq:3.6} in Theorem~A.
Together with the earlier observations \eqref{eq:3.21}--\eqref{eq:3.26} this shows Theorem~\ref{th:3} and consequently Theorem~\ref{th:2} and Theorem~\ref{th:1}, the conjecture of Erd\H{o}s, P\'olya and Tur\'an. \hfill $\blacksquare$

\noindent
J\'anos Pintz\\
R\'enyi Mathematical Institute\\
of the Hungarian Academy of Sciences\\
Budapest, Re\'altanoda u. 13--15\\
H-1053 Hungary\\
e-mail: pintz.janos@renyi.mta.hu

\end{document}